\documentclass[preprint]{elsarticle}

\usepackage{amsfonts}
\usepackage[mathscr]{eucal}
\usepackage[usenames]{color}

\newtheorem{thm}{Theorem}
\newtheorem{cor}{Corollary}
\newtheorem{lem}[thm]{Lemma}
\newtheorem{pro}[thm]{Proposition}
\newtheorem{rem}{Remark}

\newdefinition{definition}{Definition}
\newproof{pf}{Proof}
\newproof{pop}{Proof of Proposition \ref{propABC}}
\newproof{pot1}{Proof of Theorem \ref{teoFIInilpotente}}
\newproof{pot2}{Proof of Theorem \ref{Walchermejorado}}
\newproof{pot3}{Proof of Theorem \ref{Walcher:b2=0}}
\newproof{pot4}{Proof of Theorem \ref{th:apl}}
\newproof{pot5}{Proof of Theorem \ref{th:q=2}}
\newproof{poc1}{Proof of Corollary \ref{cor1}}
\newproof{poc2}{Proof of Corollary \ref{cor2}}

\newcommand{\Z}{\mathbb{Z}}

\newcommand{\Natural}{\mathbb{N}}

\newcommand{\C}{\mathbb{C}}
\newcommand{\bean}{\begin{eqnarray*}}
\newcommand{\eean}{\end{eqnarray*}}
\newcommand{\bea}{\begin{eqnarray}}
\newcommand{\eea}{\end{eqnarray}}

\def\llave#1{\left\{ #1 \right\}}

\def\fracp#1#2{{\textstyle{\frac{#1}{#2}}}}

 \def\qh{\mathscr{P}}
\def\qhf{\mathscr{F}}
 \def\QH{\mathcal{Q}}
 
\def\t{{\mathbf t}}

\def\F{{\mathbf F}}

\def\X{{\mathbf X}}
\def\D{{\mathbf D}}
\def\sn{{\mathrm{Sn}}}
\def\cs{{\mathrm{Cs}}}

\def\LD#1#2{L_{#1}#2}

 \def\x{{\mathbf x}}
 
\def\pref#1{(\ref{#1})}

\def \findemo{\hspace*{.1cm}$\Box$ \\}

\def \aiif{AIIF}

\begin{document}
{\normalfont

\title{Non-formally integrable centers admitting an algebraic inverse integrating factor}
\author{A. Algaba}
\author{N. Fuentes}
\author{C. Garc\'{\i}a}
\author{M. Reyes}

\address{Department of Mathematics. Faculty of Experimental Sciences.\\
Avda. Tres de Marzo s/n, 21071 Huelva, Spain} \ead{colume@uhu.es}


\begin{abstract}
We study the existence of a class of inverse integrating factor for a family  of non-formally integrable systems, in general,  whose lowest-degree quasi-homogeneous term is a Hamiltonian vector field. Once the existence of an inverse integrating factor is established, we characterize the systems having a center. Among others, we characterize the centers of the systems whose lowest-degree quasihomogeneous term is 
$(-y^3,x^3)^T$ with an algebraic inverse integrating factor.   
\end{abstract}


\maketitle

\section{Introduction and statement of the main results.}

One of the classic problems in the qualitative theory of the planar analytic systems is to characterize when a monodromic point (singular point which is surrounded by orbits of the system) is a center or a focus. This problem, so-called center problem, has been solved theoretically for  a nondegenerate singular point (systems whose linear part evaluated at singular point has two imaginary eigenvalues non-zero) and for the nilpotent case. Nowadays, the problem remains still unsolved for the  
remaining case, i.e. the systems with linear part identically zero at singular point, so-called  degenerate singular point.
   
One of the main tools used for characterizing the nondegenerate and nilpotent centers   has been the computation of a normal form, see Poincar\'e \cite{po1881}, Moussu \cite{mo1982}. 
It is not strange to think that a possible solution might be given by means of the theory of normal forms
 for the degenerate case.
 
Another problem related to the center problem is, once the monodromy is established,  to determine the existence of an analytic first integral. So, for instance, for a nondegenerate singular point, the analytic integrability and center problems are equivalent. 
Otherwise, the existence of a first integral is a sufficient condition but it is not necessary for the singular point to be a center. 

In this context, the existence of an integrating factor or an inverse integrating factor enable us to provide information about both center and integrability problems.

For more details about the relevance of the presence of an inverse integrating factor in a neighborhood of a singular point see \cite{gagigr2010,gagigr2011,gillvi1996} and references therein.

In this paper mainly we focus on the problem of characterizing, by means of the theory of normal forms, when a system has an inverse integrating factor in a neighborhood of the singular point. Once the existence of an inverse integrating factor and the monodromy of the origin have been established, we determine if the origin is either a center or a focus.

We consider an autonomous system
 \begin{equation} \label{Sistema}
\dot{\x}=\F(\x)=(P(\x),Q(\x))^T, \ \x \in \mathbb{C}^2,
 \end{equation}
where $\F$ is a formal planar vector field defined in a
neighborhood of the origin $U \subset \mathbb{C}^2$ having a
singular point at the origin, i.e., $\F({\bf 0})={\bf 0}$ and $P,Q\in \mathbb{C}[[x,y]]$ (algebra of the power series in $x$ and
$y$ with coefficient in $\mathbb{C}$).

A non-null $\mathcal{C}^1$ class function $V$ is an inverse
integrating factor of system \pref{Sistema} (or also of $\F$) on
$U$ if satisfies the linear partial differential equation
$\LD{\F}{V}=\mbox{div}(\F)V,$ being $\LD{\F}{V}:=P\partial
V/\partial x+Q\partial V/\partial y$, the Lie derivative of $V$
respect to $\F$, and  $\mbox{div}(\F):=\partial P/\partial
x+\partial Q/
\partial y,$ the divergence of $\F.$ This name for $V$ comes from the fact that $V^{-1}$ defines on $U\setminus
\{V=0\}$ an integrating factor of system \pref{Sistema}, i.e.
$\F/V$ is divergence-free. So, if system (\ref{Sistema}) has an formal inverse integrating factor
$V$ then it is formally integrable on  $U\setminus
\{V=0\}.$ For more details about the relation between the integrability and the
inverse integrating factor see \cite{algare2012,chgigill1999}.

We are interested in characterizing degenerate systems which
have an algebraic inverse integrating factor over
$\mathbb{C}((x,y))$ (which will be named AIIF) where
$\mathbb{C}((x,y))$ denotes the quotient field of the algebra of
the power series $\mathbb{C}[[x,y]]$. In this sense, the only results we know
are Walcher \cite{walcher2003} where is claimed its existence for non-degenerate cusp nilpotent singularity, and 
Algaba {\it et. al.}
\cite{algare2011} where is characterized all nilpotent systems  having an \aiif. 

Given $\t=(t_1,t_2)$ non-null  with $t_1$ and
$t_2$ non-negative integer numbers without common factors, 
we will denote by $\mathscr{P}^{\t}_{k}$ to the vector space of quasi-homogeneous polynomials of type $\t$ and
degree $k,$ i.e. 
$$\mathscr{P}^{\t}_{k}=\{f\in\mathbb{C}[x,y]: f(\varepsilon^{t_1}x,\varepsilon^{t_2}y)=\varepsilon^kf(x,y)\},$$
and by 
$$\mathcal{Q}^{\t}_{k}=\{  \F=(P,Q)^T :  P\in\mathscr{P}^{\t}_{k+t_1},\  Q\in\mathscr{P}^{\t}_{k+t_2}  \}$$
to the vector space of the quasi-homogeneous
polynomial vector fields of type $\t$ and degree $k.$
Any vector field can be expanded into quasi-homogeneous terms of
type $\t$ of successive degrees. Thus, the vector field $\F$
 can be written in the form $$\F=\F_r+\F_{r+1}+\cdots,$$ for some $r\in\Z,$
where $\F_j=(P_{j+t_1},Q_{j+t_2})^T\in \mathcal{Q}^{\t}_{j}$ and $\F_r\not\equiv {\bf 0}.$ Such expansions
will be expressed  as  $\F=\F_r+\mbox{q-h.h.o.t.},$ where ''q-h.h.o.t." means ''quasi-homogeneous higher order terms."

If we select the type $\t=(1,1)$, we are using in fact the Taylor expansion, but in general, each term in
the above expansion involves monomials with different degrees.

Given $h\in\mathscr{P}^{\t}_{r+|\t|},$ we define the linear operator 
\begin{eqnarray} \label{derLie} \ell_{j} &:&
\qh_{j-r}^{t} \longrightarrow \qh_j^{\t} \nonumber \\ && \mu_{j-r}
\longrightarrow \ell_{j}(\mu_{j-r}):=\frac{\partial h}{\partial x} \frac{\partial \mu_{j-r}}{\partial y}-\frac{\partial h}{\partial y} \frac{\partial \mu_{j-r}}{\partial x},
\end{eqnarray} 
(Poisson bracket of $h$ and $\mu_{j-r}$) and  denote by $\mbox{Cor}(\ell_j)$ a complementary subspace to the range of the linear operator $\ell_{j}$.

We also define $\qhf^{\t}_{r+|t|}$ as the set of all $h \in \mathscr{P}^{\t}_{r+|t|}$ satisfying:
\begin{itemize}
\item[{\bf H1}] the factorization of $h$ on
$\mathbb{C}[x,y]$ has only  simple factors,
\item[{\bf H2}]  $h\qh_{j}^{t}$ is a complementary subspace to the range of $\ell_{r+|\t|+j}$ for all $j.$
\end{itemize}
In this paper, fixed $h\in \qhf^{\t}_{r+|\t|},$   we deal with the systems of the form
\begin{eqnarray}\label{campo}
\dot{\x}=\X_h+\mbox{q-h.h.o.t.},
\end{eqnarray}
where 
$\X_{h}:=(-\partial h/\partial y , \partial
h/\partial x)^T\in\mathcal{Q}^{\t}_{r},$
i.e. a class of systems which can be
considered as perturbations of a
Hamiltonian system whose Hamiltonian function $h$ is a quasi-homogeneous function. 

This class of systems is a wide family and  contains, among others, to the  non-degenerate saddle ($h=xy$), linear center ($h=x^2+y^2$), the nilpotent systems of the form 
 $(\dot{x},\dot{y})=(y, a x^{n})+ \mbox{q-h.h.o.t}$ with $a\ne 0$ ($h=\frac{a}{n+1}x^{n+1}-\frac{y^2}{2}$).

There are two main reasons for imposing that  $h$ belongs to $\qhf_{r+|\t|}^t$. 
On the one hand, if  {\bf H1} holds, a cyclicity of the co-ranges of the operators $\ell_j$ appears.
Concretely, 
\begin{equation}\label{ciclicidad}\mbox{Cor}(\ell_{j+r+|\t|}) = h \mbox{Cor}(\ell_{j}),\ \ \mbox{
for all}\  j>r\ \mbox{with}\  \qh_{j-r}^{\t}\neq\{0\},
\end{equation}
 see \cite{algaga2009}.  Algaba {\it et. al.} \cite{alfugare2014} provide an orbital equivalent normal form up any order for the system (\ref{campo}).
This normal form is
\begin{eqnarray}\label{FNgen}
\dot{\x}=\X_h + \X_g+\mu\D_0,
\end{eqnarray}
(we have denoted $\D_0:=(t_1x,t_2y)^T\in \QH_{0}^{\t}$) 
being $g=\sum _{j\ge 1}g_{r+|\t|+j}$ with $g_{r+|\t|+j}\in \mbox{Cor}(\ell_{r+|\t|+j})\setminus h\qh_{j}^{\t}$ for $j\le r$ or $j>r$ such that $\qh_{j-r}^{\t}=\{0\},$  and  $\mu = \sum_{j> r} \mu_j, \: \mu_j\in \mbox{Cor}(\ell_j).$ 

Moreover, if $\mu_j\equiv 0,$ for all $j,$ then  system \pref{campo} is formally orbital equivalent to a
 Hamiltonian system and, in such case, it is a formally integrable system. Otherwise, from Algaba et al. \cite{algaga2009}, the system is non-formally integrable.\\
 
On the other hand, the condition {\bf H2} on $h$ implies that $g\equiv 0,$ i.e. in this paper we limit to studying the systems whose normal form is a perturbation of a Hamiltonian vector field with dissipative vector fields. \\
The following theorem summarizes the above results.
\begin{thm}[\cite{alfugare2014,algaga2009}]\label{teofn-int}
We consider system \pref{campo} with $h\in \qhf_{r+|\t|}^{\t}$. It holds that:
\begin{enumerate}
\item System \pref{campo} is formally orbital equivalent to 
$\dot{\x}=\X_h +\mu\D_0,$ with $\mu = \sum_{j> r} \mu_j$ with  $\mu_j\in \mbox{Cor}(\ell_j).$
\item System \pref{campo} is formally integrable if and only if it is formally orbital equivalent to $\dot{\x}=\X_h.$
\end{enumerate}
\end{thm}
The main result of this paper is stated in the next theorem.
\begin{thm} \label{teofiialg} System \pref{campo} with $h\in \qhf_{r+|\t|}^{\t}$ 
has an \aiif\ (algebraic inverse integrating factor over
$\mathbb{C}((x,y))$)  if and only if it is formally orbital
equivalent either to ${\dot \x}=\X_h$ (formally integrable system) or to 
\begin{equation}\label{Sisfiialg}(\dot{x},\dot{y})^T=\X_h+\mu_{r+N}\D_0,
\end{equation}
with $N$ a natural number and $\mu_{r+N}\in\mbox{Cor}(\ell_{r+N})\setminus \{0\}$ (non-formally integrable system).   
Moreover, 
the \aiif\ is $(h+\mbox{q-h.h.o.t.})^{1+N/(r+|\t|)}$,
up to a multiplicative constant. 
\end{thm}
As a consequence, it has the main result of Algaba {\it et al.} \cite{alfugare2014}.
\begin{cor} \cite[Theorem 2]{alfugare2014}
Under the conditions of Theorem \ref{teofiialg}, system \pref{campo} has a formal inverse integrating factor (it belongs to $\mathbb{C}[[x,y]]$, algebra of the power series in $x$ and
$y$ with coefficient in $\mathbb{C}$)  if and only if
it is formally orbital equivalent either to ${\dot \x}=\X_h$ (formally integrable system) or to system \pref{Sisfiialg}
with $N$ a multiple of $r+|\t|$ (non-formally integrable system).
\end{cor}
\begin{rem}\label{remunicidad}
From \cite[Theorem 3.19]{algaga2009}, system \pref{campo} is formally integrable if and only if it is formally orbital equivalents to $\dot{{\bf x}}=\X_h,$ i.e. there exist a diffeomorphism $\Phi$ and a function $\eta$ on $U\subset \C^2$ with $\mbox{det}D\Phi$ has no zero on $U$ and $\eta({\bf 0})\ne 0,$ such that 
$\Phi_*(\eta\F)=\X_h,$ where we have denoted as $\Phi_*$ to the push-forward defined by $\Phi.$ 
As $\X_h$ is a Hamiltonian vector field, $f(h)$ is a first integral for any $f$ non-constant. In particular, it is an inverse integrating factor. So,  the pull-back $\Phi^*$ brings $f(h)$ to the inverse integrating factor of system \pref{campo}, V=$f(h+\cdots)+\cdots$ i.e. it is not unique. Also, if $f(0)\ne 0,$ $V$ would be a formal inverse integrating factor with $V(0,0)\ne 0.$
\end{rem}
We study the monodromic and center problems of system \pref{campo}. For the monodromy problem, it has the following result. 
\begin{pro} \label{promono} The origin of system \pref{campo} with $h\in \qhf_{r+|\t|}^{\t}$is monodromic if and only if $h$ is only zero at the origin.
\end{pro}
 
 We note that if the origin is a monodromic point and the system is formally integrable, then the origin is a center. 
Last on, we state the result which gives title to this work where it characterizes the centers of the non-formally integrable systems \pref{campo} having an \aiif.
\begin{thm} \label{thmcenterfocus} We assume that the origin of system \pref{campo} with  $h\in \qhf_{r+|\t|}^{\t}$ is monodromic and  it is formally orbital equivalent to the non-formally integrable  system \pref{Sisfiialg}. Then, the origin
is:
\begin{enumerate}
\item a center, if $I=0,$
\item an unstable focus, if $sig(h)I>0,$
\item a stable focus, if $sig(h)I<0,$ 
\end{enumerate}
being $I=\int_{h=sig(h)}\mu_{r+N}.$
\end{thm}

\section{Some examples and applications}
In this section we show several families of systems \pref{campo} with $h\in \qhf_{r+|\t|}^{\t}$ where the origin is or not monodromic. For the non-monodromic case, we determine the systems with an \aiif.  For the monodromic case, we also characterize the centers admitting an \aiif.

In order to determine if a quasi-homogeneous function holds the condition {\bf H2}, we need to describe the sets $\mathscr{P}^{\t}_{k}$ of quasi-homogeneous polynomials according the type $\t=(t_1,t_2).$
The following result provides bases for these spaces.
\begin{lem}\label{indices}
Fixed $\t=(t_1,t_2)$, it has that:
\begin{enumerate}
\item $\mathscr{P}^{\t}_{0}=\mbox{span}\{ 1\}.$ 
\item if $t_1=1,$ for every $t_2\ge 1,$ the sets $\mathscr{P}^{\t}_{k}$ are non-trivial spaces for all $k$,
\item $\mathscr{P}^{\t}_{k}=\{0\},$ if $k\notin  \mathcal{I}^{\t},$
\item if $k>t_1t_2-|\t|,$ then $k\in \mathcal{I}^{\t},$ i.e. $\mathscr{P}^{\t}_{k}$ is a non-trivial space. 
\item $\mathscr{P}^{\t}_{k}=\mbox{span}\{ x^{k_1+t_2(k_3-j)}y^{k_2+t_1j} : j=0,\ldots,k_3\},$ if $k\in  \mathcal{I}^{\t}\setminus \{0\},$
\end{enumerate}
being $\mathcal{I}^{\t}=\{ k=k_1t_1+k_2t_2+k_3t_1t_2\in \Natural : k_1,k_2,k_3\in \Natural, k_1<t_2, k_2<t_1\}.$
\end{lem}  
 Table \ref{tab:indices} shows the sets $\Natural\setminus \mathcal{I}^{\t}$, that is, the degrees $l$ such that $\mathscr{P}^{\t}_{l}$ is a trivial set, for $t_2\le 5.$

\begin{table}[h]
\caption{Sets $\Natural\setminus \mathcal{I}^{\t}$ for $t_2\le 5.$}
\centering
\begin{tabular}{ll}
\noalign{\medskip}
\hline
$\Natural\setminus \mathcal{I}^{(1,t_2)}=\emptyset$& \\ 
$\Natural\setminus \mathcal{I}^{(2,3)}=\{ 1\} $ & $\Natural\setminus \mathcal{I}^{(2,5)}=\{ 1,3\}$\\ 
$\Natural\setminus \mathcal{I}^{(3,4)}=\{ 1,2,5\} $ & $\Natural\setminus \mathcal{I}^{(3,5)}=\{ 1,2,4,7\}$\\ 
$\Natural\setminus \mathcal{I}^{(4,5)}=\{ 1,2,3,6,7,11\} $ & \\ \hline
\end{tabular}
\label{tab:indices}
\end{table}
\begin{rem}
By \pref{ciclicidad}, the condition {\bf H2}, by assuming {\bf H1}, is equivalent to  $h\qh_{j}^{t}$ is a complementary subspace to the range of $\ell_{r+|\t|+j}$ for all $j\le r,$ or $j>r$ satisfying  $\qh_{j-r}^{t}=\{0\}$, that is, {\bf H2} holds if it satisifies a finite number of conditions.
\end{rem}
\begin{rem} \label{rem1t2}
From above lemma, the number of trivial spaces $\mathscr{P}^{\t}_{k}$ is a finite number, and by \pref{ciclicidad}, only it is enough the computation of a certain number of co-ranges, concretely, from $r+1$ to $n_0+r+|\t|-1$ (with $n_0:=1+r$ if $\Natural\setminus \mathcal{I}^{\t}$ is an empty set, or $n_0:=1+r+{\rm{max}}\{ \Natural\setminus \mathcal{I}^{\t}\},$ otherwise) for obtaining the normal form of \pref{campo}. So,  if $h\in\qhf_{r+|\t|}^{\t},$ the normal form \pref{FNgen}
provided in Algaba {\it et. al.} \cite{alfugare2014}
is 
\begin{eqnarray} \label{campofinal}\dot{\x}=\X_h+  \sum_{j=r+1}^{n_0+r+|\t|-1}\eta_j^{(0)}\D_0 + 
\sum_{i=1}^{\infty}\sum_{j=0}^{r+|\t|-1}\eta_{j+n_0}^{(i)}h^i\D_0,\end{eqnarray}
with $\eta_j^{(i)}\in\mbox{Cor}(\ell_j).$  
Moreover $r+1 \leq n_0\leq r+1+ max\{0,t_1t_2-|\t|\}$.
\end{rem}


{\sc {\bf A) Perturbations of  Hamiltonian quadratic systems.}} These systems 
 can be written as 
\begin{equation}\label{sisquad}(\dot{x},\dot{y})^T=\X_h+\mbox{q-h.h.o.t.}\ \ h=ax^3+bx^2y+cxy^2+dy^3.\end{equation}
That is, $\t=(1,1)$ and $r=1.$ 

From Proposition \ref{promono}, the origin of these systems is non-monodromic. We focus on our study in characterizing the systems \pref{sisquad} with an \aiif.
   
For $d\ne 0,$  without loss of generality,  we can assume $c=0$ and $d=1,$ the polynomial $h$ has only simple factors if $27a^2+4b^3\ne0,$ 
and by Lemma \ref{indices}, the sets $\mathscr{P}^{\t}_{j}$ are non-trivial spaces for all $j$.
Table \ref{table:quad}  shows the range and co-range of the operator $\ell_{j},\ j=2,3,4$ for system \pref{sisquad} with $d\ne 0.$ It is easy to check that $h\in \qhf_{3}^{(1,1)}.$
\begin{table}[h]
\caption{Range and co-range of operator $\ell_{j} $ for system  \pref{sisquad}. }
\label{table:quad}
\begin{tabular}{l}
\noalign{\medskip} \hline \\
Range($\ell_{2}$)=span\{$-bx^2-3y^2,3ax^2+2bxy$\}\\ 
If $a\ne 0,$ Cor($\ell_{2}$)=span\{$xy$\}. If $a=0$, Cor($\ell_{2}$)=span\{$x^2$\} \\ \noalign{\medskip} 
Range($\ell_{3}$)=span\{$-2bx^3-6xy^2,6ax^3+4bx^2y-3h,6ax^2y+4bxy^2$\}\\ 
 Cor($\ell_{3}$)=span\{$h$\} \\ \noalign{\medskip} 
 Range($\ell_{4}$)=span\{$3bx^4+9x^2y^2,-9ax^4-6bx^3y+6xh,$ \\
\hspace{3cm}$ -9ax^3y-6bx^2y^2+3yh$\}\\ 
  Cor($\ell_{4}$)=span\{$xh,yh$\}  \\ \noalign{\medskip} \hline
\end{tabular}
\end{table}

The normal form \pref{campofinal} of system \pref{sisquad} becomes
\begin{equation}\label{fnsisquad}
(\dot{x},\dot{y})^T=(-bx^2-3y^2,3ax^2+2bxy)^T+\sum_{j\ge 0}f_j(x,y,h)h^j\D_0,
\end{equation}
with $\D_0=(x,y)^T$ and $f_j\in\mbox{span}\{h,xh,yh,xyh\}$ if $a\ne 0$, or $f_j\in\mbox{span}\{h,xh,yh,x^2h\}$ if $a=0.$

Applying Theorem \ref{teofiialg}, we get the following result.
\begin{thm}
System \pref{sisquad} has an \aiif\  if and only if is formally orbital equivalent to one of the following systems:
\begin{enumerate}
\item $\dot{\x}=\X_h.$ It admits an \aiif\ of the form $g(h+\mbox{q-h.h.o.t.})$ with $g$ any nonzero function. In particular, there are inverse integrating factors nonzero at the origin. 
\item  $\dot{\x}=\X_h+\alpha_{3j}h^j\D_0,\ \alpha_{3j}\ne0, j\ge 1.$ The \aiif\ is $(h+\mbox{q-h.h.o.t.})^{j+2/3}.$
\item  $\dot{\x}=\X_h+(\alpha_{3j+1}x+\beta_{3j+1}y)h^j\D_0,\ (\alpha_{3j+1},\beta_{3j+1})\ne(0,0),  j\ge 1.$ The \aiif\ is $(h+\mbox{q-h.h.o.t.})^{1+j},$ i.e. it is a formal inverse integrating factor,
\item  $\dot{\x}=\X_h+\alpha_{3j+2}xyh^j\D_0$ if $a\ne 0$, or  $\dot{\x}=\X_h+\alpha_{3j+2}x^2h^j\D_0$ if $a=0$ with $\alpha_{3j+2}\ne0, j\ge 1.$ The \aiif\ is $(h+\mbox{q-h.h.o.t.})^{j+4/3}.$
\end{enumerate}
\end{thm}

{\sc {\bf B) Perturbations of nilpotent Hamiltonian systems.}} We consider the  nilpotent systems whose quasi-homogeneous expansion is of the form
\begin{equation} \label{Sisnilpotente}
(\dot{x},\dot{y})^T=(y,\sigma x^n)^T+\mbox{q-h.h.o.t.}\qquad \sigma=\pm 1.
\end{equation}
From Proposition \ref{promono}, the origin is not monodromic if and only if $n$ even, or $n$ odd and $\sigma=1$.

Algaba et al. \cite{algare2011} give the following result, by characterizing the systems
\pref{Sisnilpotente} which admit an \aiif.
\begin{thm} \label{teoFIInilpotente} System \pref{Sisnilpotente}
has an \aiif\ if and only if it is formally orbital
equivalent to
\begin{equation}\label{Sisnilfiia}(\dot{x},\dot{y})^T=(y,\sigma x^n)^T+\alpha_M^{(L)}x^Mh^L f(h)\D_0,
\end{equation}
with $h=2\sigma x^{n+1} -(n+1) y^2,\ \D_0=(2x,(n+1)y)^T,\
\alpha_M^{(L)}$ a real number, $f$ a function with $f(0)=1,\ L$ a
non-negative integer, and $M\in\llave{0,1,\dots,n-1}$ if $L>0$ or
$M\in\llave{\lfloor(n+1)/2\rfloor,\dots,n-1}$ if $L=0$.

Moreover, if $\alpha_M^{(L)}\neq 0$, then the system
\pref{Sisnilpotente} is not formally integrable, and if it admits an \aiif,
the \aiif\ is $(h+\mbox{q-h.h.o.t.})^{\fracp{2M+n+3}{2(n+1)}+L}$,
up to a multiplicative constant. Otherwise, if $\alpha_M^{(L)}=
0$, system \pref{Sisnilpotente} is formally integrable.
\end{thm}

If $n$ is even, system \pref{Sisnilpotente} has a formal inverse
integrating factor if and only if $\alpha_M^{(L)}=0,$ since
otherwise the number $\fracp{2M+n+3}{2(n+1)}$ is non-integer and
hence the inverse integrating factor is not formal. Therefore, as a consequence of Theorems \ref{teofn-int} and \ref{teoFIInilpotente}, it has the following result provided in Algaba et al. \cite{algare2012}.
\begin{thm} System \pref{Sisnilpotente} with $n$ even, 
has a formal inverse integrating factor if and only if it is formally integrable.
\end{thm}

Now, we analyze the center problem for  system \pref{Sisnilpotente} admitting an \aiif.
We assume that the origin is monodromic, i.e. $n$ odd ($n=2m-1$) and $\sigma=-1.$
These systems are
\begin{equation}\label{nilimpar}
(\dot{x},\dot{y})^T=(y,  -x^{2m-1})^T+\mbox{q-h.h.o.t.},\ \ m\ge 1,
\end{equation}
The first quasi-homogeneous term of the right-hand side  of \pref{nilimpar} is $\X_h\in\mathcal{Q}^{\t}_{m-1}$ with $\t=(1,m),\ h=\frac{1}{2m}x^{2m}+\frac{1}{2}my^2\in\mathscr{P}^{\t}_{2m}.$

We get the following result which characterizes the centers of the systems \pref{nilimpar} having an \aiif.
\begin{thm}\label{thm:nilcenter}
We assume that system \pref{nilimpar} has an \aiif. Then, the origin is a center if and only if it is formally orbital equivalent to a system invariant to the symmetry $(x,y,t)\rightarrow  (-x,y,-t).$
\end{thm}
\noindent {\sc Proof of Theorem \ref{thm:nilcenter}}. 
From Theorem \ref{teoFIInilpotente} 
 if system \pref{nilimpar} has an \aiif\ then it is formally orbital equivalent either to $(\dot{x},\dot{y})^T=(y,  -x^{2m-1})^T$ which is a center, or to 
\begin{equation}\label{Sisnilfiiaimpar}(\dot{x},\dot{y})^T=(y,-x^{2m-1})^T+Ax^Mh^Lf(h)\D_0,
\end{equation}
with $\D_0=(x,my)^T,\
A$ a real number non-zero, $f$ a function con $f(0)=1$, $L$ a
non-negative integer, and $M\in\llave{0,1,\dots,2m-2}$ if $L>0$ or
$M\in\llave{m,m+1,\dots,2m-2}$ if $L=0$.

By applying Theorem \ref{teofiialg}, we obtain a further reduction of the normal form 
\pref{Sisnilfiiaimpar} 
of  system \pref{nilimpar}, it which consists in assuming $f(h)$ identically one.

In order to get the centers, it is enough to compute the integral $I$ given by Theorem \ref{thmcenterfocus}. 
In this case $I=A\int_{0}^T \cs^{M}(\theta)d\theta,$ where $(\cs(\theta),Sn(\theta))^T$ is the solution of the initial value problem
$$\frac{\rm{d}{\bf x}}{\rm{d}\theta}=\X_h({\bf x}),\quad {\bf x}(0)=(1,0)^T,$$
and $T$ is a minimal period of both functions.

It is known that the integral $I$ is different from zero if and only if $M$ even. So, we arrive to a system invariant to the symmetry $(x,y,t)\rightarrow  (-x,y,-t).$

The sufficient condition is trivial. \findemo

We study the form of the \aiif's of system \pref{nilimpar}.
For $n=2m-1,$ the number $\fracp{2M+n+3}{2(n+1)}$ is
natural if $M=(2k-1)m-1$ with $k$ natural. By imposing that $M\le 2m-2,$
it has that $k=1$ and $M=m-1.$  So, we have the following result.
\begin{thm}\label{thm:centerimpar}
The origin of the system \pref{nilimpar} is a non-formally integrable center admitting a formal inverse integrating factor if and only if system \pref{Sisnilfiiaimpar} is formally orbital equivalent to 
 \begin{equation}\label{Sisnilfiiaformalcenter}(\dot{x},\dot{y})^T=(y,- x^{4k-1})^T+Ax^{2k-1}h^L\D_0,\ \ L\ge 1,\ A\ne 0.
\end{equation}
\end{thm}
Consequently, the centers of systems \pref{nilimpar} having an \aiif,  are formally orbital equivalent to time-reversible systems but no all of them have a formal inverse integrating factor.\\

{\sc {\bf  C) Quadratic nilpotent generalized systems.}} We consider the degenerate systems of the form
\begin{equation}\label{nilgen2}(\dot{x},\dot{y})^T=(y^2+\sum_{j\ge 3}P_j(x,y), \sum_{j\ge 3}Q_j(x,y))^T,
\end{equation}
with $P_j$ and $Q_j$ homogeneous polynomials of degree $j$ and $Q_3(1,0)\ne 0$ (without loss of generality, we can assume $Q_3(1,0)=1$). 
We write $P_j(x,y)=\sum_{j=m+n} a_{mn}x^my^n, \ Q_j(x,y)=\sum_{j=m+n} b_{mn}x^my^n.$
The quasi-homogeneous expansion with respect to $\t=(3,4)$ of system \pref{nilgen2} is of the form 
\begin{equation}\label{nilgen23}
(\dot{x},\dot{y})^T=(y^2, x^3)^T+\mbox{q-h.h.o.t.},
\end{equation}
i.e., system \pref{campo} for $r=5,\ h=x^4/4-y^3/3.$ 

From Proposition \ref{promono}, the origin of these systems is non-monodromic since $h$ does not preserve the sign. So, we focus on our study in characterizing the systems \pref{sisquad} with an \aiif.

Note that $h$ has only simple factors, $n_0=11$ (see Table \ref{tab:indices}).
Table \ref{table:nilgen23} shows the range and co-range of $\ell_{j}$ for $6\le j\le 22$ and $j\in\Natural\setminus \mathcal{I}^{\t}.$ 

\begin{table}[h]
\caption{Range and co-range of operator $\ell_{j} $ for system \pref{nilgen23}}
\begin{tabular}{l}
\noalign{\medskip} \hline \\
Range($\ell_{6}$)=span\{$0$\},  Cor($\ell_{6}$)=span\{$x^2$\} \\ \noalign{\medskip} 
Range($\ell_{7}$)=span\{$0$\},  Cor($\ell_{7}$)=span\{$xy$\} \\ \noalign{\medskip} 
Range($\ell_{8}$)=span\{$y^2$\},  Cor($\ell_{8}$)=span\{$0$\} \\ \noalign{\medskip} 
Range($\ell_{9}$)=span\{$x^3$\},  Cor($\ell_{9}$)=span\{$0$\} \\ \noalign{\medskip} 
Range($\ell_{10}$)=span\{$0$\},  Cor($\ell_{10}$)=span\{$x^2y$\} \\ \noalign{\medskip} 
Range($\ell_{11}$)=span\{$xy^2$\},  Cor($\ell_{11}$)=span\{$0$\} \\ \noalign{\medskip} 
Range($\ell_{12}$)=span\{$7x^4-12h$\},  Cor($\ell_{12}$)=span\{$h$\} \\ \noalign{\medskip} 
Range($\ell_{13}$)=span\{$x^3y$\}, Cor($\ell_{13}$)=\{$0$\}  \\ \noalign{\medskip} 
Range($\ell_{14}$)=span\{$x^2y^2$\}, Cor($\ell_{14}$)=\{$0$\} \\ \noalign{\medskip} 
Range($\ell_{15}$)=span\{$x^3-6xh$\}, Cor($\ell_{15}$)=span\{$xh$\}  \\ \noalign{\medskip} 
Range($\ell_{16}$)=span\{$11x^4y-12yh$\}, Cor($\ell_{16}$)=span\{$yh$\} \\ \noalign{\medskip} 
Range($\ell_{17}$)=span\{$x^3y^2$\},  Cor($\ell_{17}$)=\{$0$\}  \\ \noalign{\medskip} 
Range($\ell_{18}$)=span\{$13x^6-36x^2h$\}, Cor($\ell_{18}$)=span\{$x^2h$\} \\ \noalign{\medskip} 
Range($\ell_{19}$)=span\{$7x^5-12xyh$\}, Cor($\ell_{19}$)=span\{$xyh$\}  
\\ \noalign{\medskip} 
Range($\ell_{22}$)=span\{$17x^6y-9x^2yh$\}, Cor($\ell_{22}$)=span\{$x^2yh$\}  \\ \noalign{\medskip} \hline 
\end{tabular}
\label{table:nilgen23}
\end{table}

As above, we observe that $h\in \qhf_{12}^{(3,4)}.$
 So, the normal form \pref{campofinal} of system \pref{nilgen23} becomes
\begin{equation}\label{fnnilgen2q}
(\dot{x},\dot{y})^T=(y^2, x^3)^T+\sum_{j\ge 0}f_j(x,y,h)h^j\D_0,
\end{equation}
with $\D_0=(3x,4y)^T$ and $f_j\in\mbox{span}\{x^2,xy,x^2y,h,xh,yh\}.$

As a consequence of Theorem \ref{teofiialg}, we get the following result which characterizes the systems \pref{nilgen2} with an \aiif.
\begin{thm}\label{thm:nilgen2}
System \pref{nilgen2} has an \aiif\  if and only if, it is formally orbital equivalent to one of the following systems:
\begin{enumerate}
\item $\dot{\x}=\X_h.$ The \aiif\ is $g(h+\mbox{q-h.h.o.t.})$ with $g$ any nonzero function (in particular, there are inverse integrating factors nonzero at the origin).
\item  $\dot{\x}=\X_h+\alpha_{12j+6}x^2h^j\D_0.$ The \aiif\ is $(h+\mbox{q-h.h.o.t.})^{1+j+1/12}.$
\item  $\dot{\x}=\X_h+\alpha_{12j+7}xyh^j\D_0.$ The \aiif\ is $(h+\mbox{q-h.h.o.t.})^{1+j+1/6}.$
\item  $\dot{\x}=\X_h+\alpha_{12j+10}x^2yh^j\D_0.$ The \aiif\ is $(h+\mbox{q-h.h.o.t.})^{1+j+5/12}.$
\item  $\dot{\x}=\X_h+\alpha_{12j+12}h^{j+1}\D_0.$ The \aiif\ is $(h+\mbox{q-h.h.o.t.})^{1+j+7/12}.$
\item  $\dot{\x}=\X_h+\alpha_{12j+15}xh^{j+1}\D_0.$ The \aiif\ is $(h+\mbox{q-h.h.o.t.})^{1+j+10/12}.$
\item  $\dot{\x}=\X_h+\alpha_{12j+16}yh^{j+1}\D_0.$ The \aiif\ is $(h+\mbox{q-h.h.o.t.})^{1+j+11/12},$
\end{enumerate}
with $\alpha_k\ne 0$ and $j\ge 0.$
\end{thm}
We claim that the \aiif's of the non-formally integrable systems \pref{nilgen2} are algebraic but no formal. Consequently, we get the following result. 
\begin{pro}
System \pref{nilgen2} is formally integrable if and only if it admits a formal inverse integrating factor.
\end{pro}
Next, we give necessary conditions for the existence of an \aiif\ for system \pref{nilgen2}.
The first two coefficients of  the right-hand side of \pref{fnnilgen2q} are
\begin{eqnarray}
\alpha_6=3a_{30}+b_{21},\label{alpha6}\\
\alpha_7=13(a_{21}+b_{12})+(3a_{30}+b_{21})(4a_{30}-3b_{21}).\label{alpha7}
\end{eqnarray}
These coefficients of the quasi-homogeneous normal form have been  obtained by using the procedure given in Algaba {\it et al.}  \cite{alfrgaga2003}.\\
From Theorem \ref{thm:nilgen2}, we deduce the following result.
\begin{pro}
System \pref{nilgen2}  with $3a_{30}+b_{21}\ne 0$ is not formally integrable. Moreover,  if it has an \aiif, then $13(a_{21}+b_{12})+(3a_{30}+b_{21})(4a_{30}-3b_{21})=0$ and the \aiif\  is equal to $(4y^3-3x^4+q-h.h.o.t.)^{13/12}\mbox{exp}(u),$ for some series $u$ which is unique up to an additive constant. 
\end{pro}

We study a particular case of systems \pref{nilgen2}. We consider the family of systems \pref{nilgen2} with $P_j=Q_j\equiv 0$ for $j>3,$ and $P_3(1,0)=0$ ($a_{30}=0$), that is,
\begin{equation}\label{nilgen2esp}(\dot{x},\dot{y})^T=(y^2, x^3)^T+
(a_{21}x^2y+a_{12}xy^2+a_{03}y^3, b_{21}x^2y+b_{12}xy^2+b_{03}y^3)^T.\end{equation}
we get the following result.
\begin{pro} \label{teofiialgesp} We assume that system \pref{nilgen2esp} has an \aiif. It has that:
\begin{enumerate}
\item if $b_{21}\ne 0$, then $13(a_{21}+b_{12})=3b_{21}^2,$ (non-formally integrable case),
\item if $b_{21}=0,$ then system \pref{nilgen2esp} has a formal inverse integrating factor (integrable case).\\
Moreover, in such a case, system \pref{nilgen2esp} is one of the following systems
\begin{enumerate}
\item $b_{21}=a_{21}+b_{12}=a_{12}+3b_{03}=0,$ (Hamiltonian case).
\item  $a_{21}=a_{03}=b_{21}=b_{12}=0, a_{12}+3b_{03}\ne 0,$ (non-Hamiltonian, not axis-reversible case).
\end{enumerate}
\end{enumerate}

\end{pro}
\noindent {\sc Proof of Proposition \ref{teofiialgesp}}. 
First part follows from above proposition. We assume that $b_{21}=0$ ($\alpha_6=0$). If $a_{21}+b_{12}\ne 0$ ($\alpha_7\ne 0$),  It is easy to check that $\alpha_{10}, \alpha_{12}, \alpha_{15}$ and $\alpha_{16}$ are not zero simultaneously. Therefore, from Theorem \ref{teofiialg}, system \pref{nilgen2esp} does not have an \aiif.    
Otherwise, $a_{21}+b_{12}=0$ ($\alpha_6=\alpha_7=0$). The coefficient $\alpha_{10}$ is $\alpha_{10}=(3b_{12}-4a_{21})(a_{12}+3b_{03})$. If it is not zero, the following coefficients under the cancellation of the above ones are 
\begin{eqnarray*}
\alpha_{12}&=&(a_{12}+3b_{03})(98a_{03}+(3b_{12}-4a_{21})^2),\\
\alpha_{15}&=&(3b_{12}-4a_{21})^2(a_{12}+3b_{03})(5b_{03}-4a_{12}),\\
\alpha_{16}&=&(3b_{12}-4a_{21})(a_{12}+3b_{03})((289/1372)(3b_{12}-4a_{21})^3\\
&&\hspace{3cm}
+(11/25)(a_{12}+3b_{03})^2),\\
\alpha_{18}&=&(3b_{12}-4a_{21})^2(a_{12}+3b_{03})^3.
\end{eqnarray*}
Thus, $\alpha_{18}$ is different from zero and therefore system \pref{nilgen2esp} does not have an \aiif.\\
Otherwise, $(3b_{12}-4a_{21})(a_{12}+3b_{03})=0$ ($\alpha_6=\alpha_7=\alpha_{10}=0$). If $a_{12}+3b_{03}=0,$ system  \pref{nilgen2esp} is a Hamiltonian system whose Hamiltonian is a polynomial inverse integrating factor and a first integral. So, the system is formally integrable and it has a formal inverse integrating factor (family 2.(a)). If $a_{12}+3b_{03}\ne 0$ and $3b_{12}-4a_{21}=0$, it has that 
\begin{eqnarray*}
\alpha_{12}&=&(a_{12}+3b_{03})a_{03},\\
\alpha_{15}&=&a_{03}(a_{12}+3b_{03})(11b_{03}-8a_{12}),\\
\alpha_{16}&=&a_{03}^2(a_{12}+3b_{03}).
\end{eqnarray*}
If $a_{03}\ne 0,$ then $\alpha_{12}$ and $\alpha_{16}$ are different from zero. So, the existence of an \aiif\ arrives to $a_{03}=0,$ i.e. family 2.(b). It is straightforward to check that 
$$\begin{array}{ll}V =&1+(a_{12}+3b_{03})x+(3/2)b_{03}(a_{12}+3b_{03})x^2\\ \noalign{\medskip}
&-(1/2)b_{03}(a_{12}-3b_{03})(a_{12}+3b_{03})x^3\\ \noalign{\medskip}
&+(1/2)a_{12}b_{03}(-3b_{03}+2a_{12})(a_{12}-3b_{03})x^4\\ \noalign{\medskip}
&-(1/2)b_{03}(-b_{03}+a_{12})(a_{12}-3b_{03})(-3b_{03}+2a_{12})y^3\\ \noalign{\medskip} 
&-(1/2)b_{03}(-b_{03}+a_{12})(a_{12}-3b_{03})(-3b_{03}+2a_{12})a_{12}xy^3,\end{array}$$
is a polynomial inverse integrating factor for family 2.(b) which is 1 at origin. 
Thus, $$H=-\int{P/V \: dy} + \int{\left(Q/V + \frac{\partial}{\partial
x}\int{P/V\:dy}\right)dx}$$ is a formal first integral
defined in a neighborhood of the origin. Therefore, the system is formally integrable.
\findemo

\begin{rem}
If there exists an \aiif\ of system \pref{campo} which does not have the form $(h+\mbox{q-h.h.o.t.})^{1+j/(r+|\t|)}$,  up to a multiplicative constant, for a certain $j$, then the system is formally integrable.
For instance, $V=(y^3/3-x^4/4-3\lambda x^5)^{6/5}$ is an inverse integrating factor of 
$({\dot x},{\dot y})^T=(y^2+60\lambda xy^2,x^3+100\lambda y^3)^T,$ and from Proposition \ref{teofiialgesp}, it is an formally integrable system. 
\end{rem}

Last on, we study the problem for the systems \pref{nilgen2} given by
\begin{equation}\label{nilgen2espa30}(\dot{x},\dot{y})^T=(y^2, x^3)^T+
(a_{30}x^3, b_{21}x^2y+b_{03}y^3)^T,\end{equation} with $a_{30}\ne 0$, (case $a_{30}=0,$ studied before). 
It has the following result.
\begin{pro} \label{teofiialgespa30} System \pref{nilgen2espa30}, with $a_{30}\ne 0,$  has an \aiif\ if and only if it satisfies:
\begin{enumerate}
\item $3a_{30}+b_{21}=b_{03}=0$, (Hamiltonian system), or 
\item  $3b_{21}-4a_{30}=0$ and $b_{03}=0,$ (non-formally integrable system).\\
Moreover, in this case, the \aiif\ is $V=(4y^3-3x^4)^{13/12}.$
\end{enumerate}
\end{pro}
\noindent {\sc Proof of Proposition \ref{teofiialgespa30}}. 
We assume that system \pref{nilgen2espa30} with $a_{30}\ne 0,$  has an \aiif. 
The first two coefficients of the quasi-homogeneous normal form of  \pref{nilgen2espa30} are
given by \pref{alpha6} and \pref{alpha7} for $a_{21}=b_{12}=0.$ 
Therefore, if $3a_{30}+b_{21}\ne 0,$ it arrives to $4a_{30}-3b_{21}=0.$ In such case, the following coefficient of the normal form is $\alpha_{10}=a_{30}^2b_{03}.$ So, $b_{03}=0.$
It is easy to check that $V=(4y^3-3x^4)^{13/12}$ is an \aiif\ of the system.\\
Otherwise, $3a_{30}+b_{21}=0.$ In this case, $\alpha_6$ and $\alpha_7$ are zero and $\alpha_{10}=a_{30}^2b_{03}.$ This arrives to $b_{03}=0,$ i.e. it is a Hamiltonian system.
\findemo

{\sc {\bf D) Systems of the form $(-y^3,x^3)^T+\mbox{q-h.h.o.t.}$}}
The systems are
\begin{equation}\label{y3x3}
\dot{\x}=\X_h+\mbox{q-h.h.o.t.}
\end{equation}
with $h=x^4/4+y^4/4,\ \t=(1,1)$ and $r=2.$
From Lemma \ref{indices}, the sets $\mathscr{P}^{\t}_{j}$ are non-trivial spaces for all $j$,
hence $n_0=1+r=3$. So, in order to get a normal form, it is enough to compute the sets Cor($\ell_j$),\ $j=3,4,5,6,$ which are given in Table \ref{table:y3x3}.
\begin{table}[h]
\caption{Range and co-range of operator $\ell_{j} $ for system  \pref{y3x3}. }
\label{table:y3x3}
\begin{tabular}{l}
\noalign{\medskip} \hline \\
Range($\ell_{3}$)=span\{$x^3,y^3$\}\\ 
Cor($\ell_{3}$)=span\{$x^2y,xy^2$\}. \\ \noalign{\medskip} 
Range($\ell_{4}$)=span\{$xy^3,x^4+2h,x^3y$\}\\ 
 Cor($\ell_{4}$)=span\{$x^2y^2,h$\} \\ \noalign{\medskip} 
 Range($\ell_{5}$)=span\{$x^2y^3,3x^5+8xh,3x^4y+4yh,x^3y^2$\} \\
  Cor($\ell_{5}$)=span\{$xh,yh$\} \\ \noalign{\medskip} 
 Range($\ell_{6}$)=span\{$x^3y^3,x^6+3x^2h,x^5y+2xyh, x^4y^2+y^2h$\} \\
  Cor($\ell_{6}$)=span\{$x^2h,xyh,y^2h$\}  \\ \noalign{\medskip} \hline
\end{tabular}
\end{table}

We note that $h\in\qhf_{2}^{(1,1)}.$ A normal form of system \pref{y3x3} is
\begin{equation}\label{fny3x3}
(\dot{x},\dot{y})^T=(-y^3,x^3)^T+\sum_{j\ge 0}f_j(x,y,h)h^j\D_0,
\end{equation}
with $\D_0=(x,y)^T$ and $f_j\in\mbox{span}\{x^2y,xy^2,h,x^2y^2,xh,yh,x^2h,xyh,y^2h\}$.

 From Proposition \ref{promono}, the origin is a monodromic singular point. 
In order to characterize the centers of system \pref{fny3x3}, it is necessary to compute the value of the integrals $I_{n,k}=\int_{0}^T \cs^n(\theta)Sn^k(\theta) d\theta,\ n,k\in\{ 0,1,2\},$ being 
$g(\theta)=(\cs(\theta),\sn(\theta)),\ \theta\in[0,T)$ a parameterization of the closed curve $h=1,$ where $(\cs(\theta),\sn(\theta))^T$ is the solution of the initial value problem
$$\left\{ \begin{array}{l}
\frac{\rm{d}\cs\theta}{\rm{d}\theta}=-\sn^3\theta,\\ \noalign{\medskip}
\frac{\rm{d}\sn\theta}{\rm{d}\theta}=\cs^3\theta,\end{array}\right.$$
with $(\cs(0),\sn(0))=(1,0),$
and $T$ is a minimal period of both functions.

We cite some properties of these integrals. For the shake of shortness, we prefer to avoid its proof in this paper.
\begin{lem}\label{lemintegrales} For every $n,k\ge 0,$ it holds:
\begin{enumerate}
\item $I_{2n+1,k}=I_{n,2k+1}=0,$
\item $I_{2n+2,2k+2}=\frac{(2n+1)(2k+1)}{4(n+k+2)(n+k+1)}I_{2n,2k}.$
\end{enumerate}
\end{lem}
So, 
$$I_{1,0}=I_{0,1}=I_{1,1}=I_{2,1}=I_{1,2}=0,\quad I_{2,0}=I_{0,2},\quad I_{0,0}=8I_{2,2}.$$

Applying Theorems \ref{teofiialg} and \ref{thmcenterfocus},
 we have the following result.
\begin{thm}\label{thm:y3x3}
System \pref{y3x3} has an \aiif\  if and only if, it is formally orbital equivalent to 
\begin{enumerate}
\item $\dot{\x}=\X_h.$ The \aiif\ is $g(h)$ with $g$ any nonzero function (in particular, there are inverse integrating factors nonzero at the origin).\\
In this case, the origin is a center. 
\item  $\dot{\x}=\X_h+(\alpha_{4j+3}x^2y+\beta_{4j+3}x^2y)h^j\D_0,\ (\alpha_{4j+3},\beta_{4j+3})\ne(0,0), j\ge 0.$ The \aiif\ is $(h+\mbox{q-h.h.o.t.})^{2+j+1/4}.$\\
In this case, the origin is a center. 
\item  $\dot{\x}=\X_h+(\alpha_{4j+4}h+\beta_{4j+4}x^2y^2)h^j\D_0,\ (\alpha_{4j+4},\beta_{4j+4})\ne(0,0), j\ge 0.$ The \aiif\ is $(h+\mbox{q-h.h.o.t.})^{2+j+1/2}.$\\
In this case, the origin is a center if and only if $8\alpha_{4j+4}+\beta_{4j+4}=0$.
 \item  $\dot{\x}=\X_h+(\alpha_{4j+5}xh+\beta_{4j+5}yh)h^{j}\D_0,\ (\alpha_{4j+5},\beta_{4j+5})\ne(0,0), j\ge 0.$ The \aiif\ is $(h+\mbox{q-h.h.o.t.})^{2+j+3/4}.$\\
In this case, the origin is a center. 
\item  $\dot{\x}=\X_h+(\alpha_{4j+6}x^2h+\beta_{4j+6}xyh+\gamma_{4j+6}y^2h)h^{j}\D_0,\ (\alpha_{4j+6},\beta_{4j+6},\gamma_{4j+6})\ne(0,0,0), j\ge 0.$ The \aiif\ is $(h+\mbox{q-h.h.o.t.})^{3+j},$ i.e. it is formal.\\
In this case, the origin is a center if and only if $\alpha_{4j+6}+\gamma_{4j+6}=0.$ 
\end{enumerate}
\end{thm}

We analize the system 
\begin{equation}\label{y3x3moussu}
(\dot{x},\dot{y})=(y^3, -x^3+c_3x^2y^2+c_4xy^3).
\end{equation}
This system for $c_3=1/2$ and $c_4=0$ has been studied in Moussu \cite{mo1982} by showing that it is a degenerate analytic center without formal first integral and 
Gin\'e and Peralta-Salas \cite{gipe-sa2012} have proved that the system does not admit a formal inverse integrating factor. 

The first coefficients of the normal form \pref{fny3x3} are $\alpha_3=2c_3,\ \beta_3=3c_4.$
If both $c_3$ and $c_4$ are zero, the system is a Hamiltonian system whose first integral is $h=x^4+y^4.$ Otherwise, the system is not formally integrable. Moreover,
\begin{itemize}
\item if $c_3.c_4\ne0,$ the coefficients  of fourth order of the normal form are $\alpha_4=0$ and $\beta_4=2c_3c_4.$ Thus, from Theorem \ref{thm:y3x3}, it does not have an \aiif,
 \item if $c_3=0$ and $c_4\ne0,$ the coefficients of fourth order are zero but $\alpha_5=6/5c_4^3.$ And if $c_3\ne0$ and $c_4=0,$ it has that $\beta_5=16/45c_3^3.$ Therefore, from Theorem \ref{thm:y3x3}, it does not have an \aiif. 
\end{itemize} 
Summarizing, 
\begin{thm} System \pref{y3x3moussu} with $(c_3,c_4)\ne (0,0)$ does not admit an algebraic inverse integrating factor.
\end{thm}

\section{Proofs of the main results.}\label{sec:proofs}

The following result we will be used for the proof of Theorem \ref{teofiialg}  is an adjustment of Proposition 10 and Proposition 13 of \cite{algare2011}:
\begin{lem}
\label{lemfiialg} Let system $\dot{\x}=\X_h+\mu \D_0,$ where  the factorization of
$h\in\mathscr{P}_{r+|\t|}^\t$ on $\mathbb{C}[x,y]$ only has simple
factors and $\mu=\sum_{j\ge N}\mu_{r+j}$ with
 $\mu_{r+j}\in \mbox{Cor}(\ell_{j}),$ for all $j\ge N>0$ and $\mu_{r+N}\not\equiv 0$.
 If $V$ is an \aiif\  of the system, then
$V=(\sum_{j\ge 1}b_jh^j)^{(r+N+|\t|)/(r+|\t|)},$ with $b_1=1$, is
the unique \aiif\  up to a multiplicative constant.\\ 
Moreover, the real numbers $b_j$ verify the recursive relation
\begin{eqnarray}\label{ecuFII} 0=\sum_{i=0}^{j-1}
\left[\fracp{N+(1+i)(r+|\t|)}{r+|\t|+N}-(j-i)\right]b_{j-i}h^{j-i}\mu_{r+N+i(r+|\t|)}.\end{eqnarray}
 Furthermore, if $\mu=\lambda
f(h)+\nu$ with $\lambda\in\mbox{Cor}(\ell_{r+N})\setminus\llave{0}$,
 $f$ a scalar function, $f(0)=1,$ and $\nu=\sum_{j>N}\nu_{r+j}$,
$\nu_{r+j}\in\mbox{Cor}(\ell_{j})$, $\nu_{r+N+l(r+|\t|)}\equiv 0$
for all non-negative integer $l$, then under these conditions,
the system has an \aiif\  if and only if $\nu\equiv 0$.
\end{lem}

\noindent {\sc Proof of Theorem \ref{teofiialg}}. 
We prove the necessity.  From Theorem \ref{teofn-int}, a normal form of system \pref{campo} is of the form $\dot{\x}=\X_h + \mu\D_0,$ with 
 $\mu = \sum_{j> r} \mu_j, \: \mu_j\in \mbox{Cor}(\ell_j).$

If $\mu_j\equiv 0,$ for all $j$ then  system \pref{campo} is formally orbital equivalent to a
 Hamiltonian system and, in such case, it is proved that  system \pref{campo} has a formal inverse integrating factor.\\  

Otherwise, let $N=\mbox{min}\{ j,\ \mu_{r+j}\ne 0\}.$
By \cite[Theorem 13]{alfugare2014}, system (\ref{FNgen})  is formally orbital equivalent to 
$\dot{\x}=\X_h+ \mu_{r+N}\D_0+\sum_{j>
N}\tilde{\mu}_{r+j}\D_0,$  with $\tilde{\mu}_{j}\in \mbox{Cor}(\ell_{j}^{(2)}),$ a complementary subspace to the range of the linear operator $\ell_{k}^{(2)}\: :\: \qh_{k-r}^{t} \times \mbox{Ker}(\ell_{k-N})\longrightarrow \qh_k^{\t}
$
defined by
$$\begin{array}{ll}
\ell_{k}^{(2)}(\mu_{k-r}, \alpha h^{l_1}):=\ell_{k}(\mu_{k-r})+ \alpha \mu_{r+N}
h^{l_1}, &\mbox{ if }\: l_2=0,\\ \noalign{\medskip}
\ell_{k}^{(2)}(\mu_{k-r},0):=\ell_{k}(\mu_{k-r}), & \mbox{ if }\: l_2\ne 0.
\end{array}$$
being
$k=r+N+l_1(r+|\t|)+l_2,$ with $0\le l_2<r+|\t|.$ 

From Lemma \ref{lemfiialg},  $V=(h+\sum_{j>1}b_jh^j)^{(r+N+|\t|)/(r+|\t|)}$ with $b_j$ verifying (\ref{ecuFII}). We see that $b_j=0$ and $\tilde{\mu}_{r+N+(j-1)(r+|\t|)}=0,$ for any $j>1.$
Indeed, we assume the contrary and let $j_0=min\{ j>1:\ b_j\ne 0\}.$ By  (\ref{ecuFII}), for $j=j_0$, it has that
$$ b_{j_0}h^{j_0}\mu_{r+N}-\fracp{r+|\t|}{r+|\t|+N}h\tilde{\mu}_{r+N+(j_0-1)(r+|\t|)}=0.$$
Consequently,
$b_{j_0}h^{j_0-1}\mu_{r+N}\in\mbox{Cor}(\ell_{r+N+(j_0-1)(r+|\t|)}^{(2)})\setminus\{0\},$
but also
$$b_{j_0}h^{j_0-1}\mu_{r+N}=\ell_{r+N+(j_0-1)(r+|\t|)}^{(2)}(0,b_{j_0}h^{j_0-1}),$$
which is a contradiction. So, $b_{j_0}=0$ and $\tilde{\mu}_{r+N+(j-1)(r+|\t|)}=0,$ for all $j>1.$

Applying Lemma \ref{lemfiialg}, for $\lambda=\mu_{r+N},\ f(h)=1$ and 
$\nu_{r+j}=\mu_{r+j},$ it has that $\mu_{r+j}=0,$ for any $j>N.$

We prove that the condition is sufficient. If $\mu_{r+N}\equiv 0,$ the polynomial
$h^m$ with $m$ any natural number, is a
polynomial  first integral and, in particular, it is an inverse
integrating factor. Thus, if
we perform the transformation which brings $\dot{{\bf x}}=\X_h$  to the system
\pref{campo}, then system \pref{campo} admits an
\aiif (in fact, it is formal) but  it is not unique modulus a multiplicative constant.\\
In the case, $\mu_{r+N}\not\equiv 0,$ we see that $V(h)=h^{\fracp{r+|\t|+N}{r+|\t|}}$ is an \aiif\ of system \pref{Sisfiialg}. Indeed, 
applying Euler theorem for quasi-homogeneous function, i.e.
$\LD{\D_0}{f}=sf$ with $f\in\mathscr{P}_{s}^{\t},$ it has that the Lie derivative of $V$ respect to $\F=\X_h+\mu_{r+N}\D_0$ is 
\bean
\LD{\F}{V}&=&V'(h)\LD{\F}{h}=(r+|\t|)\mu_{r+N} hV'(h)=(r+|\t|+N)\mu_{r+N} h^{\fracp{r+|\t|+N}{r+|\t|}}\eean and
\bean
\mbox{div}(\F)&=& \mbox{div}(\mu_{r+N} \D_0)= \LD{\D_0}{\mu_{r+N}}+|\t|\mu_{r+N}=
(r+|\t|+N)\mu_{r+N}.
\eean
So, $\LD{\F}{V}-V\mbox{div}(\F)=0,$ that is, $V$ is  an \aiif\ of system \pref{Sisfiialg} (formal if $N$ is a multiple of $r+|\t|$). Thus, 
the system \pref{campo} has the \aiif,
$(h+\mbox{q-h.h.o.t.})^{\fracp{r+|\t|+N}{r+|\t|}}$, up to a
multiplicative constant.
\findemo

\noindent {\sc Proof of Proposition \ref{promono}}. 
As $h\in\qhf_{r+|\t|}^{t},$ we can write in a compact form 
$h=c\prod_{j=1}^nf_{j}\prod_{j=1}^mg_j,$ where $f_j=x,\ y$ or $y^{t_1}-\lambda_jx^{t_2},\ j=1,\dots, n$, 
$g_j(x,y)=(y^{t_1}-a_jx^{t_2})^2+ b_j^2x^{2t_2},\ j=1,\dots, m$ with $c,\lambda_j, a_j$ and $b_j$ real numbers and $\lambda_j, b_j$ non-zero, for all $j$.

We see the necessary condition. We assume the contrary one. Thus, 
$f_j$  is a factor of $h$. From Proposition 8 of \cite{algare2011b}, there exists a orbit of the system which leaves or enters at origin. Consequently, the origin is not monodromic.

On the other hand, the sufficient condition follows from Proposition 6 of \cite{algare2011b}. \findemo
\\

\noindent {\sc Proof of Theorem \ref{thmcenterfocus}}. 
We assume that system \pref{campo} is formally orbital equivalent to system \pref{Sisfiialg}.

We also can assume that $h(x,y)$ is positive for all $(x,y)\ne(0,0)$ since from Proposition \ref{promono} if the origin is monodromic $h$ preserves its sign and if $h$ is non-positive by changing the time $t$ by $-t$, $h$ becomes $-h.$  

A parameterization of the closed curve $h=1$ is $g(\theta)=(\cs(\theta),\sn(\theta)),\ \theta\in[0,T)$ where $(\cs(\theta),\sn(\theta))^T$ is the solution of the initial value problem
$$\frac{\rm{d}{\bf x}}{\rm{d}\theta}=\X_h({\bf x}),\quad {\bf x}(0)=(1,0)^T,$$
and $T$ is a minimal period of both functions.

We consider the transformation 
\begin{equation}\label{coorgeneralizada} 
x=u^{t_1}\cs(\theta),\qquad y=u^{t_2}\sn(\theta),
\end{equation}
where $u>0, \theta\in[0,T).$\\
 Differentiating \pref{coorgeneralizada} with respect to time, we get ${\bf \dot{x}}= \frac{1}{u}\D_0{\bf \dot{u}}+\frac{1}{u^r}\X_h\dot{\theta}.$ From this, we obtain
\begin{equation} \label{cambioaux}
{\bf \dot{x}}\wedge \X_h=\frac{1}{u}\D_0\wedge \X_h \dot{u},\qquad 
\D_0\wedge {\bf \dot{x}}=\frac{1}{u^r}\D_0\wedge \X_h \dot{\theta}.
\end{equation}
 We note that $\D_0\wedge\X_h(x,y)=\nabla h(x,y)\cdot \D_0=(r+|\t|)h(x,y)\ne 0,\ $ for all $(x,y)\ne (0,0).$\\
For system \pref{Sisfiialg}  it has ${\bf \dot{x}}\wedge \X_h=u^{r+N}\mu_{r+N}(\cs(\theta),\sn(\theta))\D_0\wedge\X_h$ and  $\D_0\wedge {\bf \dot{x}}=\frac{1}{u^r}\D_0\wedge \X_h.$ So, system \pref{Sisfiialg} is
\begin{equation}\label{siscorpol}
\left\{\begin{array}{l} \dot{u}=u^{r+N+1}\mu_{r+N}(\cs(\theta),\sn(\theta)),\\ \dot{\theta}=u^r. \end{array}\right .
\end{equation}
It can be further simplified by rescaling the time by $\rm{d}t=\frac{1}{u^r}\rm{d}\tau,$ which yields 
\begin{equation}\label{siscorpolred}
\left\{\begin{array}{l} u' =\frac{\rm{d}u}{\rm{d}\tau}=u^{N+1}\mu_{r+N}(\cs(\theta),\sn(\theta)),\\ \noalign{\medskip}
 \theta'=\frac{\rm{d}\theta}{\rm{d}\tau}=1. \end{array}\right .
\end{equation}
Finally, the change $z=-\frac{1}{N}u^N$ transforms the system into 
\begin{equation}\label{siscorpolred2}
\left\{\begin{array}{l} z' =\mu_{r+N}(\cs(\theta),\sn(\theta)),\\ \noalign{\medskip}
 \theta'=1. \end{array}\right .
\end{equation}
Poincar\'e map for system \pref{siscorpolred2} is
$\Pi(z_0)=z(T,z_0)=z_0+I.$ From this, the result follows. \findemo




{\bf Acknowledgments.} This work has been partially supported by
{\it Ministerio de Ciencia y Tecnolog\'{\i}a, Plan Nacional I+D+I}
co-financed with FEDER funds, in the frame of the project
 MTM2010-20907-C02-02, and by {\it
Consejer\'{\i}a de Educaci\'on y Ciencia de la Junta de
Andaluc\'{\i}a} (FQM-276 and P12-FQM-1658).


\begin{thebibliography}{00}
%



\bibitem{alfrgaga2003} {\sc A. Algaba; E. Freire;
E. Gamero; C. Garc\'{\i}a,} {\em Quasihomogeneous normal forms,} J.
Comput. Appl. Math, \textbf{150}, (2003), 193-216.



\bibitem{alfugare2014} {\sc A. Algaba; N. Fuentes; C. Garc\'{\i}a; M. Reyes,} {\em A class of non-integrable systems admitting an inverse integrating factor,} Journal Math. Anal. App. \textbf{420}, (2014), 2, 1439-1454.

\bibitem{algaga2009}{\sc A. Algaba; E. Gamero; C. Garc\'\i a,}
{\em The integrability problem for a class of planar systems,} Nonlinearity \textbf{22}, (2009), 2, 395-420.




\bibitem{algare2011}{\sc A. Algaba; C. Garc\'\i a; M. Reyes}, {\em Nilpotent systems admitting an algebraic  inverse integrating factor over $\C((x,y))$,} Qualitative Theory of Dynamical Systems, \textbf{10}, 2, 303-316, (2011).

\bibitem{algare2011b}{\sc A. Algaba; C. Garc\'\i a; M. Reyes}, {\em Characterization of a monodromic singular point of a planar vector field}, Nonlinear Analysis, \textbf{74}, 5402-5414, (2011).

\bibitem{algare2012}{\sc A. Algaba; C. Garc\'\i a; M. Reyes}, {\em Existence of an inverse integrating factor, center problem and integrability of a class
of nilpotent systems,} Chaos Solitons \& Fractals, \textbf{45}, 869-878, (2012).











\bibitem{chgigill1999}
{\sc J. Chavarriga; H. Giacomini; J. Gin\'e; J. Llibre,} {\em On
the integrability of two-dimensional flows,} J. Differential
Equations \textbf{157}, (1999), 163-182.








\bibitem{gagigr2010}{\sc I. Garc\'\i a; H. Giacomini; M. Grau,} {\em The inverse integrating factor and the Poincar\'e map.} Trans. Amer. Math. Soc.
\textbf{362} (2010), 7, 3591-3612.

\bibitem{gagigr2011}{\sc I. Garc\'\i a; H. Giacomini; M. Grau,} {\em Generalized Hopf Bifurcation for planar vector fields via the inverse integrating factor.} J. Dyn. Differ. Equat., \textbf{23}, 2,(2011), 251-281.







\bibitem{gillvi1996}
{\sc H. Giacomini; J. Llibre; M. Viano,} {\em On the nonexistence,
existence and uniqueness of limit cycles,} Nonlinearity
\textbf{9}, (1996), 501-516.



\bibitem{gipe-sa2012}
{\sc J. Gin\'e; D. Peralta-Salas,} {\em Existence of inverse integrating factors and Lie symmetries for degenerate planar centers,} J. Differential Equations
\textbf{252}, (2012), 344-357.



\bibitem{mo1982} {\sc R. Moussu,} {\em Sym\'etrie et forme normaledes centres et foyers d\'eg\'en\'er\'es,} Ergodic Theory Dynam. Sys.\textbf{2} (1982), 241-251.




\bibitem{po1881} {\sc H. Poincar\'e, } {\em M\'emoire sur les courbes d\'efinies par les \'equations diff\'erentielles,} J. Math.\textbf{37}, (1881), 375-422.










\bibitem{walcher2003}
{\sc S. Walcher,} {\em Local integrating factors,} J. Lie Theory
\textbf{13}, (2003), 279-289.

\end{thebibliography}
\end{document}